\RequirePackage{fix-cm}
\documentclass[smallextended]{svjour3}       
\smartqed  
\usepackage{latexsym}
\usepackage{amsmath}
\usepackage{setspace} 
\usepackage{comment}
\usepackage{bm}
\usepackage{amsfonts}
\usepackage{graphicx}
\usepackage{graphics}
\usepackage{colortbl}
\usepackage{mathrsfs}
\usepackage{mathtools}
\usepackage{multirow}
\newcommand{\R}{\mathbb{R}}
\newcommand{\N}{\mathbb{N}}
\newcommand{\C}{\mathbb{C}}
\renewcommand{\S}{\mathbb{S}}
\renewcommand{\L}{\mathbb{L}}

\newcommand{\maxi}{\mathop{\rm maximize}}
\newcommand{\Mini}{\mathop{\rm Minimize}}

\newcommand{\subj}{\mbox{subject to}}

\newtheorem{algorithm}{Algorithm}

\begin{document}

\title{An equivalent nonlinear optimization model with triangular low-rank factorization for semidefinite programs}

\titlerunning{A new equivalent nonlinear optimization model for SDP}        

\author{
Yuya Yamakawa
\and
Tetsuya Ikegami
\and
Ellen H. Fukuda
\and
Nobuo Yamashita
}


\institute{
Yuya Yamakawa
\at
Department of Applied Mathematics and Physics, Graduate School of Informatics, Kyoto University, Yoshidahommachi, Sakyo-ku, Kyoto-shi, Kyoto 606-8501, Japan, 
\\
\email{yuya@i.kyoto-u.ac.jp}
\and
Tetsuya Ikegami
\at
Data and Science Solutions Group, Technology Group, Yahoo Japan Corporation, Kioi Tower, Tokyo Garden Terrace Kioicho, 1-3 Kioicho, Chiyoda-ku, Tokyo 102-8282, Japan, 
\\
\email{teikegami@yahoo-corp.jp}
\and
Ellen H. Fukuda
\at
Department of Applied Mathematics and Physics, Graduate School of Informatics, Kyoto University, Yoshidahommachi, Sakyo-ku, Kyoto-shi, Kyoto 606-8501, Japan, 
\\
\email{ellen@i.kyoto-u.ac.jp}
\and
Nobuo Yamashita
\at
Department of Applied Mathematics and Physics, Graduate School of Informatics, Kyoto University, Yoshidahommachi, Sakyo-ku, Kyoto-shi, Kyoto 606-8501, Japan, 
\\
\email{nobuo@i.kyoto-u.ac.jp}
}

\date{Received: date / Accepted: date}

\maketitle

\begin{abstract}
In this paper, we propose a new nonlinear optimization model to solve semidefinite optimization problems (SDPs), providing some properties related to local optimal solutions. The proposed model is based on another nonlinear optimization model given by Burer and Monteiro (2003), but it has several nice properties not seen in the existing one. Firstly, the decision variable of the proposed model is a triangular low-rank matrix, and hence the dimension of its decision variable space is smaller. Secondly, the existence of a strict local optimum of the proposed model is guaranteed under some conditions, whereas the existing model has no strict local optimum. In other words, it is difficult to construct solution methods equipped with fast convergence using the existing model. Some numerical results are also presented to examine the efficiency of the proposed model.

\end{abstract}

\section{Introduction}
In this paper, we consider the following {\it semidefinite optimization problem} (SDP):
\begin{eqnarray*} \label{sdp}
(\mbox{SDP}) \quad
\begin{array}{llll}
\displaystyle \Mini_{X \in \S^{n}} & \displaystyle \langle C, X \rangle &
\\
\rm{subject ~ to}  & \displaystyle {\cal A}(X) = b, ~ X \succeq O,
\end{array}
\end{eqnarray*}
where $\S^{n}$ denotes the set of $n \times n$ real symmetric matrices, and the operator ${\cal A} \colon \S^{n} \to \R^{m}$ is defined by ${\cal A}(X) := [ \langle A_{1}, X \rangle \cdots \langle A_{m}, X \rangle ]^{\top}$, and the matrices $C, ~ A_{1}, \ldots, A_{m} \in \S^{n}$ and the vector $b \in \R^{m}$ are given. For two matrices $P$ and $Q$ included in $\R^{p \times q}$, the inner product of them is defined by $\langle P, Q \rangle := {\rm tr}(P^{\top}Q)$, where ${\rm tr}(M)$ represents the trace of a square matrix $M$, and the superscript $\top$ indicates the transposition of a matrix or a vector. Let $\S_{++}^{n} ~ (\S_{+}^{n})$ be the set of $n \times n$ real symmetric positive (semi)definite matrices. For a matrix $M \in \S^{n}$, $M \succeq O$ and $M \succ O$ mean that $M \in \S_{+}^{n}$ and $M \in \S_{++}^{n}$, respectively.
\par
SDPs include some classes of optimization problems, such as linear programs, quadratic programs, and second-order cone programs. Moreover, they have a wide range of application fields, such as control theory, graph theory, structural optimization, combinational optimization, and so forth \cite{To01,Va96,WoSaVa00}. Until now, a lot of solution methods for SDPs have been proposed by many researchers \cite{MoOrSv14,PoReWi06,SuToYa16,To01,WeGoYi09,YaSuTo15,ZhSuTo10}. Specifically, the primal-dual interior-point method is known as one of the most popular ones, and there exist several efficient software packages in which they are implemented, such as SeDuMi \cite{St99}, SDPT3 \cite{ToToTu99}, and SDPA \cite{YaFuKo03}. They are based on the Newton method, and can solve small- or medium-scale problems very accurately. However, it may not be applicable to problems whose scale is too large. Moreover, it is well known that the Newton equation used in the primal-dual interior-point method becomes unstable in the neighborhoods of solutions, and hence it is difficult to obtain solutions with high accuracy.
\par
To solve large-scale SDPs whose $m$ (the number of equality constraints) is small and $n$ (the matrix size) is large, Burer and Monteiro have proposed an equivalent nonlinear optimization model \cite{BuMo03}. In their proposal, the decision variable $X$ is replaced with a matrix product $RR^{\top}$, where $R \in \R^{n \times r}$. Since $RR^{\top}$ is positive semidefinite, the semidefinite constraint $X \succeq O$ can be ignored, and hence the model is expressed as a usual nonlinear optimization problem. It is known that if $r \geq \max \{ r \in \N \colon \frac{r(r+1)}{2} \leq m\}$, then their model is equivalent to (SDP). Note that depending on the value $r$, the dimension of such model's decision variable space can be smaller than $\frac{n(n+1)}{2}$, which is the dimension of the original (SDP). However, the model has no strict local minimum. Therefore, the second-order sufficient conditions do not hold, and it is difficult to construct solution methods equipped with fast convergence. Moreover, when $m$ is not small enough, the dimension $nr$ is larger than $\frac{n(n+1)}{2}$.
\par
In this paper, we present a new nonlinear optimization model, which overcomes the above drawbacks of the Burer and Monteiro's model, and show that the proposed model has several nice properties associated with local optima. A remarkable point of the proposed model is that its decision variable space is the set of $n \times r$ real matrices whose upper triangular part is all zero, that is, its dimension is equal to $nr - \frac{r(r-1)}{2}$. It is smaller than $nr$, which is the dimension of the variables of the existing model. Moreover, the dimension $nr - \frac{r(r-1)}{2}$ in the proposed model is at most $\frac{n(n+1)}{2}$, that is to say, it never exceeds the dimension of the original (SDP). We also show that the existence of a strict local optimum is guaranteed under some appropriate conditions. From this fact, it is expected that second-order methods, such as sequential quadratic programming (SQP) and interior-point methods, have fast convergence to solutions. Herein, we also provide an SQP method to solve the proposed model.
\par
This paper is organized as follows. In Section $2$, we introduce some important concepts regarding (SDP) and the existing model proposed by Burer and Monteiro. In Section $3$, we propose a new nonlinear optimization model and give some properties associated with its local optima. In Section $4$, we provide an SQP method for the proposed model. Section $5$ reports some numerical experiments to confirm the efficiency of the proposed model. Finally, we make some concluding remarks in Section $6$.
\par
Throughout this paper, we use the following notation. The identity matrix and the all-ones vector are represented by $I$ and $e$, respectively, where their dimensions are defined by each context. For a vector $w \in \R^{p}$, $[w]_{i}$ denotes the $i$-th element of $w$, and $\Vert w \Vert$ is the Euclidean norm of $w$ defined by $\Vert w \Vert := \sqrt{\left\langle w, w \right\rangle} \, (\, = \sqrt{w^{\top}w})$. Let $W \in \R^{p \times q}$. We express the $(i,j)$-entry of $W$ by $[W]_{ij}$. Moreover, we write $\Vert W \Vert_{{\rm F}}$ and $\Vert W \Vert_{2}$ for the Frobenius norm and the operator norm of $W$, respectively, that is, $\Vert W \Vert_{{\rm F}} := \sqrt{\left\langle W, W \right\rangle} \, (\, = \sqrt{{\rm tr}(W^{\top}W)})$ and $\Vert W \Vert_{2} := \sup \{ \Vert Wx \Vert \colon \Vert x \Vert = 1 \}$, where ${\rm tr}(M)$ denotes the trace of a square matrix $M$. For real numbers $r_{1}, \ldots, r_{d} \in \R$ and a vector $v \in \R^{d}$, we use the notation below:
\begin{eqnarray*}
{\rm diag}(r_{1}, \ldots, r_{d}) := \left[
\begin{array}{ccc}
r_{1} & & O
\\
& \ddots &
\\
O & & r_{d}
\end{array}
\right], \quad {\rm diag}(v) := \left[
\begin{array}{ccc}
[v]_{1} & & O
\\
& \ddots &
\\
O & & [v]_{d}
\end{array}
\right].
\end{eqnarray*}
Let $U \in \S^{d}$ be a matrix. The minimum and the maximum eigenvalues of $U$ are denoted by $\lambda_{\min}(U)$ and $\lambda_{\max}(U)$, respectively. Let $\Phi$ be a mapping from $P_{1} \times P_{2}$ to $P_{3}$, where $P_{1}$ and $P_{2}$ are open sets. We express the Fr\'echet derivative of $\Phi$ as $\nabla \Phi$. Moreover, we denote the Fr\'echet derivative of $\Phi$ with respect to a variable $Z \in P_{1}$ as $\nabla_{Z} \Phi$. For a positive integer $k \in \N$, we define
\begin{eqnarray}
r_{k} := \max \left\{ r \in \N \colon \frac{r(r+1)}{2} \leq k \right\}. \label{def:rl}
\end{eqnarray}

\section{Preliminaries} \label{sec_pre}
We give some important concepts related to (SDP). Next, we introduce the existing nonlinear model for (SDP) and provide some of its properties.

\subsection{Basic facts related to SDP}
As it is well known, the dual of (SDP) can be written as
\begin{eqnarray*} \label{dualsdp}
(\mbox{DSDP}) \quad
\begin{array}{lll}
\displaystyle \maxi_{(y,Z) \in \R^{m} \times \S^{n}} &  \langle b, y \rangle
\\
\mbox{subject to}  & \displaystyle {\cal A}^{\ast}(y) + Z = C, \quad Z \succeq O,
\end{array}
\end{eqnarray*}
where ${\cal A}^{\ast} \colon \R^{m} \to \S^{n}$ is the adjoint operator of ${\cal A}$, which is defined by ${\cal A}^{\ast}(v) := \sum_{j=1}^{m} [v]_{j}A_{j}$ for all $v \in \R^{m}$. Throughout this paper, we assume the existence of $( X^{\ast}, y^{\ast}, Z^{\ast} ) \in \S^{n} \times \R^{m} \times \S^{n}$ such that
\begin{eqnarray}
{\cal A}(X^{\ast}) = b, ~ {\cal A}^{\ast}(y^{\ast}) + Z^{\ast} = C, ~ \langle X^{\ast}, Z^{\ast} \rangle = 0, ~ X^{\ast} \succeq O, ~ Z^{\ast} \succeq O, \label{KKTconditions}
\end{eqnarray}
which are called Karush-Kuhn-Tacker (KKT) conditions of (SDP). These conditions are necessary and sufficient for optimality. The following result shows the existence of a solution with a particular limited rank.

\begin{theorem} \label{equiv_rank} 
There exists an optimal solution $X^{\ast} \in \S^{n}$ of {\rm (SDP)} such that ${\rm rank}(X^{\ast}) \leq r_{m}$.
\end{theorem}
\noindent
{\it Proof.}
It follows from definition \eqref{def:rl} and \cite[Theorem 1]{BuMo03} (see also \cite{Ba95} and \cite{Pa98}).

\subsection{An existing nonlinear optimization model for SDP}
In \cite{BuMo03}, Burer and Monteiro proposed the following low-rank SDP, which has a rank constraint on the decision variable $X \in \S^{n}$:
\begin{eqnarray*} \label{lrsdp}
(\mbox{LRSDP}_{r}) \quad
\begin{array}{ll}
\displaystyle \Mini_{X \in \S^{n}} & \langle C, X \rangle
\\
\displaystyle \mbox{subject to}  & {\cal A}(X) = b, ~ X \succeq 0, ~ {\rm rank}(X) \leq r.
\end{array}
\end{eqnarray*}
Theorem~\ref{equiv_rank} ensures that if $r \geq r_{m}$, then (LRSDP$_{r}$) is equivalent to (SDP). Since an arbitrary semidefinite matrix $X \in \S^{n}$ can be rewritten as $X = RR^{\top}$ for some $R \in \R^{n \times r}$, (LRSDP$_{r}$) can be reformulated as follows:
\begin{eqnarray*} \label{nsdp}
(\mbox{NSDP}_{r}) \quad
\begin{array}{lll}
\displaystyle \Mini_{R \in \R^{n \times r}} & \displaystyle \langle C, RR^{T} \rangle
\\
\mbox{subject to}  & \displaystyle {\cal A}(RR^{\top}) = b.
\end{array}
\end{eqnarray*}
For (NSDP$_{r}$), the Lagrange function is defined as
\begin{eqnarray*}
L(R, v) := \langle C, RR^{\top} \rangle - \langle v, {\cal A}(RR^{\top}) - b \rangle = \langle C - {\cal A}^{\ast}(v), RR^{\top} \rangle + \langle b, v \rangle,
\end{eqnarray*}
where $v \in \R^{m}$ is the Lagrange multiplier. We say that a feasible point $R^{\ast}$ is stationary of (NSDP$_{r}$) if there exists $v^{\ast} \in \R^{m}$ such that
\begin{eqnarray*}
\nabla_{R} L(R^{\ast}, v^{\ast}) = 2( C - {\cal A}^{\ast}(v^{\ast}) )R^{\ast} = 0.
\end{eqnarray*}
\par
Problem (NSDP$_{r}$) has several remarkable properties. In particular, we do not have to deal directly with the semidefinite constraint, and if $r \in \N$ is small, then the number of variables decreases considerably compared with (SDP). However, since the purpose of this paper is to obtain a solution of the original (SDP), we need to clarify the relation between the global (or local) optimal solutions of (SDP) and (NSDP$_{r}$). In fact, some of these relations can be seen in the results below.

\begin{proposition} \label{opt} {\rm \cite[Proposition~2.3~and~Theorem~3.4]{BuMo05}}
Suppose that $R^{\ast}$ is a local minimum of {\rm (NSDP$_{r}$)} with $r \geq r_{m+1}$. If {\rm (SDP)} has a unique optimal solution, then it is given by $R^{\ast} (R^{\ast})^{\top}$.
\end{proposition}

\begin{proposition} \label{opt2} {\rm \cite[Proposition~3]{BuMo03}}
Suppose that $R^{\ast}$ is a stationary point of {\rm (NSDP$_{r}$)}, i.e., there exists $v^{\ast} \in \R^{m}$ such that $\nabla_{R} L(R^{\ast}, v^{\ast}) = 0$. If $C - {\cal A}^{\ast}(v^{\ast})$ is positive semidefinite, then $R^{\ast}(R^{\ast})^{\top}$ and $(C - {\cal A}^{\ast}(v^{\ast}), v^{\ast})$ are optimal solutions for {\rm (SDP)} and {\rm (DSDP)}, respectively.
\end{proposition}

\begin{remark} \label{existing_local_opt}
Let $R$ be a local optimum of (NSDP$_{r}$). Note that $RQ(RQ)^{\top} = RR^{\top}$ for an arbitrary orthogonal matrix $Q$. Then, we easily see that $RQ$ is also another local optimum of (NSDP$_{r}$). Moreover, we can select the matrix $Q$ so that $RQ$ arbitrarily approaches to $R$, and hence (NSDP$_{r}$) has no strict local optimum \cite{BuMo03}. This fact shows that it is difficult to construct fast convergent methods which solve (NSDP$_{r}$).
\end{remark}

\section{A new nonlinear optimization model for SDP} \label{sec_TNSDP}
Firstly, we propose a new nonlinear optimization model for (SDP). Secondly, we provide relations between solutions of the proposed model and (SDP), and some important properties related to local optima.
\par
To begin with, we denote by $\L^{n \times r}$ the set of lower triangular matrices in $\R^{n \times r}$, i.e.,
\begin{eqnarray*}
\L^{n \times r} := \{ S \in \R^{n \times r} \colon [S]_{ij} = 0 ~ \mbox{if} ~ i<j \}.
\end{eqnarray*}
The result below shows that a symmetric positive semidefinite matrix in $\S^{n}$ with rank smaller than $n$ has at least one triangular low-rank factorization.

\begin{proposition} \label{low_rank_factorization}
For any symmetric positive semidefinite matrix $X \in \S^{n}$ satisfying ${\rm rank}(X) = r < n$, there exists a matrix $S \in \L^{n\times r}$ such that $X = SS^{\top}$.
\end{proposition}
\noindent
{\it Proof.} Since X is a symmetric positive semidefinite matrix satisfying ${\rm rank}(X) = r < n$, there exists a matrix $U \in \R^{n\times r}$ such that $X = UU^{\top}$. Let us write the matrix $U$ as follows:
\begin{eqnarray*}
U = \left[
\begin{array}{cc}
U_{1}
\\
U_{2}
\end{array}
\right], \quad U_{1} \in \R^{r \times r}, \quad U_{2} \in \R^{(n-r) \times r}.
\end{eqnarray*}
From the QR factorization of the matrix $U_{1}^{\top}$, there exist an orthogonal matrix $Q \in \R^{r \times r}$ and an upper triangular matrix $R \in \R^{r \times r}$ such that $U_{1}^{\top} = QR$. Now, let $S$ be defined by
\begin{eqnarray*}
S := \left[
\begin{array}{c}
R^{\top}
\\
U_{2} Q
\end{array}
\right].
\end{eqnarray*}
Note that $S$ is included in $\L^{n \times r}$ because the matrix $R^{\top}$ is lower triangular. Hence, we have
\begin{eqnarray*}
X &=& \left[
\begin{array}{c}
R^{\top}Q^{\top}
\\
U_2
\end{array}
\right] \left[
\begin{array}{cc}
QR & U_2^{\top}
\end{array}
\right] 
\\
&=& \left[
\begin{array}{cc}
R^{\top} R & R^{\top} Q^{\top} U_{2}^{\top}
\\
U_{2} QR & U_{2} U_{2}^{\top}
\end{array}
\right] 
\\
&=& \left[
\begin{array}{c}
R^{\top}
\\
U_{2} Q
\end{array}
\right] \left[
\begin{array}{cc}
R & Q^{\top}U_{2}^{\top}
\end{array}
\right] 
\\
&=& SS^{\top}.
\end{eqnarray*}
Therefore, the assertion is proven. \qed
\bigskip

The above proposition guarantees that (LRSDP$_{r}$) can be rewritten as the following new nonlinear optimization model:
\begin{eqnarray*} \label{nsdp_new}
(\mbox{T-NSDP}_{r}) \quad
\begin{array}{lll}
\displaystyle \Mini_{S \in \L^{n \times r}} & \displaystyle \langle C, SS^{T} \rangle
\\
\mbox{subject to}  & \displaystyle {\cal A}(SS^{\top}) = b.
\end{array}
\end{eqnarray*}
The Lagrange function for (T-NSDP$_{r}$) is defined by
\begin{eqnarray*}
\mathcal{L}(S, w) := \langle C, SS^{\top} \rangle - \langle w, {\cal A}(SS^{\top}) - b \rangle = \langle C - {\cal A}^{\ast}(w), SS^{\top} \rangle + \langle b, w \rangle,
\end{eqnarray*}
where $w \in \R^{m}$ is the Lagrange multiplier.
\par
Similarly to the existing model (NSDP$_{r}$), we now verify how the local (global) optimal solutions of (SDP) and (T-NSDP$_{r}$) are related to each other. The following propositions can be proven in the same way to \cite{BuMo03,BuMo05}.

\begin{proposition} \label{opt_new_model} {\rm \cite[Proposition~2.3~and~Theorem~3.4]{BuMo05}}
Suppose that $S^{\ast}$ is a local minimum of {\rm (T-NSDP$_{r}$)} with $r \geq r_{m+1}$. If {\rm (SDP)} has a unique optimal solution, then it is given by $S^{\ast} (S^{\ast})^{\top}$.
\end{proposition}

\begin{proposition} \label{opt2_new_model} {\rm \cite[Proposition~3]{BuMo03}}
Suppose that $S^{\ast}$ is a stationary point of {\rm (T-NSDP$_{r}$)}, i.e., there exists $w^{\ast} \in \R^m$ such that $\nabla_{S} \mathcal{L}(S^{\ast}, w^{\ast}) = 0$. If $C - {\cal A}^{\ast}(w^{\ast})$ is positive semidefinite, then $S^{\ast}(S^{\ast})^{\top}$ and $(C - {\cal A}^{\ast}(w^{\ast}), w^{\ast})$ are optimal solutions for {\rm (SDP)} and {\rm (DSDP)}, respectively.
\end{proposition}

The proposed model (T-NSDP$_{r}$) has several advantages over the existing model (NSDP$_{r}$). One of them is that the number of variables in (T-NSDP$_{r}$) (more precisely, $nr - \frac{r(r-1)}{2}$) can be smaller than (NSDP$_{r}$)'s (i.e., $nr$). In particular, if $r = n$, then the number of variables in (NSDP$_{r}$) is larger compared to the original (SDP). Meanwhile, (T-NSDP$_{r}$) has at most $\frac{n(n+1)}{2}$ variables, and hence it is always less than or equal to (SDP)'s.
\par
The reduction of variables mentioned above also brings another nice property. As stated in Remark~\ref{existing_local_opt}, (NSDP$_{r}$) has no strict local optimum. On the other hand, (T-NSDP$_{r}$) has a strict local optimum under some appropriate conditions because it decreases the degree of freedom of variables compared with (NSDP$_{r}$). To show this fact, we recall the following well-known result.

\begin{proposition} \label{cholesky} {\rm \cite[Fact~8.9.38]{Be09}}
If $X \in \S^{N}$ is a positive definite matrix, then there exists a unique matrix $V \in \L^{N \times N}$ such that $X = VV^{\top}$ and $[V]_{jj} > 0$ for all $j \in \{ 1, \ldots, N \}$.
\end{proposition}
\noindent
Notice that this result basically shows that the Cholesky decomposition is unique. By exploiting such fact, we show the next lemma.
\begin{lemma} \label{uniqueness_lemma}
Suppose that $P \in \L^{\ell \times \ell}$ is an arbitrary full-rank matrix. Then, there exists $\delta > 0$ such that if $Q \in \L^{\ell \times \ell}$ satisfies $PP^{\top} = QQ^{\top}$, then $P = Q$ or $\Vert P - Q \Vert_{{\rm F}} \geq \delta$.
\end{lemma}
\noindent
{\it Proof.}
We show the assertion by contradiction. Let $F := PP^{\top}$. Because ${\rm rank}(P) = \ell$, it is clear that $F \succ O$, that is, $\lambda_{j}(F) > 0$ for all $j \in \{ 1, \ldots, \ell \}$. Hence, we define $\delta := \sqrt{\lambda_{\min}(F)} > 0$. Since the assertion is not true, there exists $Q_{\delta} \in \L^{\ell \times \ell}$ such that
\begin{eqnarray} \label{assumeQ}
F = Q_{\delta}Q_{\delta}^{\top}, \quad P \not = Q_{\delta}, \quad \Vert P - Q_{\delta} \Vert_{{\rm F}} < \delta .
\end{eqnarray}
Note that $[P]_{jj} \not = 0$ and $[Q_{\delta}]_{jj} \not = 0$ for all $j \in \{ 1, \ldots, \ell \}$ because $P \in \L^{\ell \times \ell}$, $Q_{\delta} \in \L^{\ell \times \ell}$, and $\ell = {\rm rank}(F) = {\rm rank}(P) = {\rm rank}(Q_{\delta})$. Now, we define diagonal matrices $G \in \R^{\ell \times \ell}$ and $H_{\delta} \in \R^{\ell \times \ell}$ satisfying
\begin{eqnarray} \label{def:GH}
[G]_{jj} := \left\{
\begin{array}{rl}
 1 & \mbox{if} ~ [P]_{jj} > 0,
\\
-1 & \mbox{if} ~ [P]_{jj} < 0,
\end{array}
\right. \quad [H_{\delta}]_{jj} := \left\{
\begin{array}{rl}
 1 & \mbox{if} ~ [Q_{\delta}]_{jj} > 0,
\\
-1 & \mbox{if} ~ [Q_{\delta}]_{jj} < 0,
\end{array}
\right. \quad j = 1, \ldots, \ell.
\end{eqnarray}
Then, there exists $\widehat{P} \in \L^{\ell \times \ell}$ such that 
\begin{eqnarray}
P = \widehat{P}G, \quad [\widehat{P}]_{jj} > 0, ~ \forall j \in \{ 1, \ldots, \ell \}. \label{mat:widehatP}
\end{eqnarray}
Similarly, there exists $\widehat{Q}_{\delta} \in \L^{\ell \times \ell}$ such that $Q_{\delta} = \widehat{Q}_{\delta}H_{\delta}$ and $[\widehat{Q}_{\delta}]_{jj} > 0$ for all $j \in \{ 1, \ldots, \ell \}$. It follows from these results that $F = \widehat{P} \widehat{P}^{\top} = \widehat{Q}_{\delta} \widehat{Q}_{\delta}^{\top}$. However, Proposition~\ref{cholesky} ensures that $\widehat{P} = \widehat{Q}_{\delta}$. As a result, we obtain 
\begin{eqnarray}
P = \widehat{P}G, \quad Q_{\delta} = \widehat{P}H_{\delta}. \label{eq:matPQ}
\end{eqnarray}
Combining \eqref{assumeQ} and \eqref{eq:matPQ} yields
\begin{eqnarray} 
\lambda_{\min}(F) = \delta^{2}
&>& \Vert P - Q_{\delta} \Vert_{{\rm F}}^{2} \nonumber
\\
&=& \Vert \widehat{P} ( G - H_{\delta} ) \Vert_{{\rm F}}^{2} \nonumber
\\
&=& {\rm tr}( \widehat{P}^{\top}\widehat{P} (G - H_{\delta})^{2} ) \nonumber
\\
&\geq& \lambda_{\min}(F){\rm tr}((G-H_{\delta})^{2}), \label{ineq:lambdaF} 
\end{eqnarray}
where the third equality follows from the fact that ${\rm tr}(AB) = {\rm tr}(BA)$ for any matrices $A$ and $B$, and the last inequality is true because ${\rm tr}(AB) \geq \lambda_{\min}(A){\rm tr}(B)$ for $A \in \S^{\ell}$ and $B \in \S_{+}^{\ell}$ \cite[Theorem 8.4.13]{Be09}. Moreover, since $\widehat{P} \in \L^{\ell \times \ell}$ is nonsingular from \eqref{mat:widehatP}, the results \eqref{assumeQ} and \eqref{eq:matPQ} mean that $O \not = P - Q_{\delta} = \widehat{P}(G-H_{\delta})$, i.e., $G \not = H_{\delta}$. Exploiting \eqref{def:GH} and \eqref{ineq:lambdaF} implies $\lambda_{\min}(F)  > \lambda_{\min}(F){\rm tr}((G-H_{\delta})^{2}) \geq 4\lambda_{\min}(F)$, that is, $0 \geq \lambda_{\min}(F)$. Therefore, this contradicts $F \succ O$. \qed
\bigskip

In the following, we provide sufficient conditions under which (T-NSDP$_{r}$) has a strict local optimum.

\begin{theorem} \label{strict_local}
Assume that {\rm (SDP)} has a unique optimal solution $X^{\ast}$ satisfying ${\rm rank}(X^{\ast}) = \ell \in [1, r]$. Suppose also that $X^{\ast}$ has the following structure:
\begin{eqnarray*}
X^{\ast} = \left[
\begin{array}{cc}
X_{1}^{\ast} & O
\\
O & O
\end{array}
\right], \quad X_{1}^{\ast} \in \S^{\ell}, \quad {\rm rank}(X_{1}^{\ast}) = \ell.
\end{eqnarray*}
If either of the following two statements holds, then $S^{\ast}$ is a strict local optimum of {\rm (T-NSDP$_{r}$)}:
\begin{itemize}
\item[(i)] $S^{\ast} \in \L^{n \times r}$ is a local optimum of {\rm (T-NSDP$_{r}$)} with $r \geq r_{m+1}$;
\item[(ii)] there exists $(S^{\ast}, w^{\ast}) \in \L^{n \times r} \times \R^{m}$ such that $\nabla_{S} \mathcal{L}(S^{\ast}, w^{\ast}) = 0$ and $C - {\cal A}(w^{\ast}) \succeq O$.
\end{itemize}
\end{theorem}
\noindent
{\it Proof.}
Firstly, we consider the case where statement~(i) holds. Proposition~\ref{opt_new_model} implies that $X^{\ast} = S^{\ast}(S^{\ast})^{\top}$ is a unique optimum of (SDP). Let $S^{\ast} \in \L^{n \times r}$ be denoted by
\begin{eqnarray*} \label{structure_Sstar0}
S^{\ast} = \left\{
\begin{array}{lll}
\left[
\begin{array}{cc}
S_{1}^{\ast} & O
\\
S_{2}^{\ast} & S_{3}^{\ast}
\end{array}
\right],  && \mbox{if} ~ 1 \leq \ell < r,
\vspace{1mm}
\\
\left[
\begin{array}{c}
S_{1}^{\ast}
\\
S_{2}^{\ast}
\end{array}
\right],  && \mbox{if} ~ \ell = r,
\end{array}
\right.
\end{eqnarray*}
where $S_{1}^{\ast} \in \L^{\ell \times \ell}$. Since $X^{\ast} = S^{\ast}(S^{\ast})^{\top}$ holds, we get
\begin{eqnarray*}
\left[
\begin{array}{cc}
X_{1}^{\ast} & O
\\
O & O
\end{array}
\right] = \left\{
\begin{array}{lll}
\left[
\begin{array}{cc}
S_{1}^{\ast}(S_{1}^{\ast})^{\top} & S_{1}^{\ast}(S_{2}^{\ast})^{\top}
\\
S_{2}^{\ast}(S_{1}^{\ast})^{\top} & S_{2}^{\ast}(S_{2}^{\ast})^{\top} + S_{3}^{\ast}(S_{3}^{\ast})^{\top}
\end{array}
\right],  && \mbox{if} ~ 1 \leq \ell < r,
\vspace{1mm}
\\
\left[
\begin{array}{cc}
S_{1}^{\ast}(S_{1}^{\ast})^{\top} & S_{1}^{\ast}(S_{2}^{\ast})^{\top}
\\
S_{2}^{\ast}(S_{1}^{\ast})^{\top} & S_{2}^{\ast}(S_{2}^{\ast})^{\top}
\end{array}
\right],  && \mbox{if} ~ \ell = r.
\end{array}
\right.
\end{eqnarray*}
Thus, we easily see that
\begin{eqnarray*}
0 = \left \{ 
\begin{array}{lll}
{\rm tr}(S_{2}^{\ast}(S_{2}^{\ast})^{\top} + S_{3}^{\ast}(S_{3}^{\ast})^{\top}) = \Vert S_{2}^{\ast} \Vert_{{\rm F}}^{2} + \Vert S_{3}^{\ast} \Vert_{{\rm F}}^{2}, && \mbox{if} ~ 1 \leq \ell < r,
\\
{\rm tr}(S_{2}^{\ast}(S_{2}^{\ast})^{\top}) = \Vert S_{2}^{\ast} \Vert_{{\rm F}}^{2}, && \mbox{if} ~ \ell = r.
\end{array}
\right.
\end{eqnarray*}
As a result, we obtain $S_{1}^{\ast} \in \L^{\ell \times \ell}$, $X_{1}^{\ast} = S_{1}^{\ast} (S_{1}^{\ast})^{\top}$, ${\rm rank}(S_{1}^{\ast}) = {\rm rank}(X_{1}^{\ast}) = \ell$, and
\begin{eqnarray} \label{structure_Sstar}
S^{\ast} = \left\{
\begin{array}{lll}
\left[
\begin{array}{cc}
S_{1}^{\ast} & O
\\
O & O
\end{array}
\right], && \mbox{if} ~ 1 \leq \ell < r,
\vspace{1mm}
\\
\left[
\begin{array}{c}
S_{1}^{\ast}
\\
O
\end{array}
\right], && \mbox{if} ~ \ell = r.
\end{array}
\right.
\end{eqnarray}
Now, it follows from Lemma~\ref{uniqueness_lemma} and ${\rm rank}(S_{1}^{\ast}) = \ell$ that there exists $\delta > 0$ such that
\begin{eqnarray} \label{delta_cond}
Q \in \L^{\ell \times \ell}, ~ QQ^{\top} = S_{1}^{\ast}(S_{1}^{\ast})^{\top} \quad \Longrightarrow \quad Q = S_{1}^{\ast} ~~ \mbox{or} ~~ \Vert Q - S_{1}^{\ast} \Vert_{{\rm F}} \geq \delta.
\end{eqnarray}
Let $S \in \L^{n \times r}$ be an arbitrary matrix satisfying 
\begin{eqnarray} \label{condition_S}
{\cal A}(SS^{\top}) = b, \quad \Vert S - S^{\ast} \Vert_{{\rm F}} < \delta, \quad S \not = S^{\ast}.
\end{eqnarray}
Note that $SS^{\top} \not= S^{\ast}(S^{\ast})^{\top}$ is a sufficient condition under which $S^{\ast}$ is a strict local optimum. Indeed, if $SS^{\top} \not= S^{\ast}(S^{\ast})^{\top}$ holds, then the uniqueness of $X^{\ast} = S^{\ast}(S^{\ast})^{\top}$ implies that $\langle C, SS^{\top} \rangle > \langle C, X^{\ast} \rangle = \langle C, S^{\ast}(S^{\ast})^{\top} \rangle$. Hence, we show $SS^{\top} \not= S^{\ast}(S^{\ast})^{\top}$ by contradiction. In the following, we consider the case where $1 \leq \ell < r$. Concerning the case where $\ell = r$, we can prove $SS^{\top} \not= S^{\ast}(S^{\ast})^{\top}$ in a similar way, and hence we omit its proof.
\par
Let $S$ be represented as follows:
\begin{eqnarray} \label{structure_S}
S = \left[
\begin{array}{cc}
S_{1} & O
\\
S_{2} & S_{3}
\end{array}
\right],  \quad S_{1} \in \L^{\ell \times \ell}, \quad S_{2} \in \R^{(n - \ell) \times \ell}, \quad S_{3} \in \R^{(n - \ell) \times (r - \ell)}.
\end{eqnarray}
Combining (\ref{structure_Sstar}), (\ref{structure_S}), and the assumption $SS^{\top} = S^{\ast}(S^{\ast})^{\top}$ yields
\begin{eqnarray} \label{block_mat}
\left[
\begin{array}{cc}
S_{1}S_{1}^{\top} & S_{1}S_{2}^{\top}
\\
S_{2}S_{1}^{\top} & S_{2}S_{2}^{\top} + S_{3}S_{3}^{\top}
\end{array}
\right] = SS^{\top} = S^{\ast}(S^{\ast})^{\top} = \left[
\begin{array}{cc}
S_{1}^{\ast} (S_{1}^{\ast})^{\top} & O
\\
O & O
\end{array}
\right].
\end{eqnarray}
Notice that $S_{1} \in \L^{\ell \times \ell}$ from (\ref{structure_S}), and that $S_{1}S_{1}^{\top} = S_{1}^{\ast}(S_{1}^{\ast})^{\top}$ from (\ref{block_mat}). It then follows from (\ref{delta_cond}) that $S_{1} = S_{1}^{\ast}$ or $\Vert S_{1} - S_{1}^{\ast} \Vert_{{\rm F}} \geq \delta$. Since condition (\ref{condition_S}) leads to $\Vert S_{1} - S_{1}^{\ast} \Vert_{{\rm F}} \leq \Vert S - S^{\ast} \Vert_{{\rm F}} < \delta$, we get $S_{1} = S_{1}^{\ast}$. Now, recall that $S \not = S^{\ast}$ by condition (\ref{condition_S}). Then, (\ref{structure_Sstar}) and (\ref{structure_S}) yield $S_{2} \not = O$ or $S_{3} \not = O$. However, we have from (\ref{block_mat}) that $\Vert S_{2} \Vert_{{\rm F}}^{2} + \Vert S_{3} \Vert_{{\rm F}}^{2} = {\rm tr}(S_{2}S_{2}^{\top} + S_{3}S_{3}^{\top}) = 0$, i.e., $S_{2} = O$ and $S_{3} = O$. Therefore, we see that $SS^{\top} \not= S^{\ast}(S^{\ast})^{\top}$.
\par
Secondly, we assume that statement~(ii) holds. Proposition~\ref{opt2_new_model} ensures that $X^{\ast} = S^{\ast}(S^{\ast})^{\top}$ is a unique solution of (SDP). Thus, we can use the same arguments from the case~(i) and this completes the proof.
\qed
\bigskip

In Theorem~\ref{strict_local}, we assume that the solution $X^{\ast}$ has a certain block structure. Although one may consider that such a structure is not generally satisfied, it can be assumed without loss of generality. In what follows, we explain this fact.
\par
We consider the case where $X^{\ast}$ does not have the structure given in Theorem~\ref{strict_local}. Since $X^{\ast} \in \S^{n}$ and ${\rm rank}(X^{\ast}) = \ell$, there exist an orthogonal matrix $U \in \R^{n \times n}$ and a diagonal matrix $D \in \S^{n}$ such that $X^{\ast} = U D U^{\top}$ and ${\rm rank}(D) = \ell$. Moreover, the diagonal matrix $D$ has the following structure:
\begin{eqnarray} \label{structure_D}
D = 
\begin{array}{rl}
\left[
\begin{array}{cc}
D_{1} & O
\\
O & O
\end{array}
\right], \quad D_{1} \in \S^{\ell}, \quad {\rm rank}(D_{1}) = \ell.
\end{array}
\end{eqnarray}
Note that $\langle C, X \rangle = \langle U^{\top}CU, U^{\top}XU \rangle$, $\langle A_{j}, X \rangle = \langle U^{\top}A_{j}U, U^{\top}XU \rangle ~(j=1,\ldots, m)$, and $X \succeq O$ if and only if $U^{\top}XU \succeq O$. Therefore, we can reformulate (SDP) as follows:
\begin{eqnarray} \label{reformSDP}
\begin{array}{llll}
\displaystyle \Mini_{\widehat{X} \in \S^{n}} & \displaystyle \langle \widehat{C}, \widehat{X} \rangle &
\\
\rm{subject ~ to}  & \displaystyle \langle \widehat{A}_{j}, \widehat{X} \rangle = [b]_{j} ~(j=1,\ldots,m), ~ \widehat{X} \succeq O,
\end{array}
\end{eqnarray}
where $\widehat{C} := U^{\top}CU$ and $\widehat{A}_{j} := U^{\top} A_{j} U ~ (j=1,\ldots,m)$. It is clear that (SDP) is equivalent to problem~(\ref{reformSDP}), and hence $U^{\top}X^{\ast}U = D$ is an optimal solution of problem~(\ref{reformSDP}). If we replace (SDP) with problem~(\ref{reformSDP}), then the solution $D$ satisfies the block structure required in Theorem~\ref{strict_local} because it satisfies (\ref{structure_D}).

\begin{remark}
A result that corresponds to Theorem~\ref{strict_local} was not considered in \cite{BuMo03,BuMo05} because the existing model (NSDP$_{r}$) has no strict local optimum as described in Remark~\ref{existing_local_opt}. By Theorem~\ref{strict_local}, it is expected that solutions can be obtained more rapidly if we apply second-order methods to (T-NSDP$_{r}$), such as SQP or interior point methods.
\end{remark}

\section{An SQP method for the nonlinear optimization models}
Let us first mention that the nonlinear optimization models (NSDP$_{r}$) and (T-NSDP$_{r}$) can be expressed as a certain quadratic equality constrained quadratic program (QECQP). We later provide an SQP method to solve QECQP based on this fact.
\par
We now show that the nonlinear models (NSDP$_{r}$) and (T-NSDP$_{r}$) can be recast into the following optimization problem:
\begin{eqnarray*} \label{reform_nonlinear_model}
\mbox{(QECQP)} \quad
\begin{array}{ll}
\displaystyle \Mini_{x \in \R^{d}} & \displaystyle f(x) := \frac{1}{2} \langle H x, x \rangle
\vspace{1mm}
\\
\mbox{subject to} & \displaystyle g_{j}(x) := \frac{1}{2} \langle G_{j} x, x \rangle - \frac{1}{2} [b]_{j} = 0 ~ (j=1,\ldots,m),
\end{array}
\end{eqnarray*}
where matrices $H$ and $G_{j} ~ (j=1,\ldots,m)$ have the following block structures:
\begin{eqnarray*}
&& H = \left[
\begin{array}{ccc}
\widetilde{H}_{1} & & O
\\
 & \ddots &
\\
O & & \widetilde{H}_{r}
\end{array}
\right], \quad \widetilde{H}_{k} \in \S^{s_{k}} ~ (k = 1,\ldots, r),
\\
&& G_{j} = \left[
\begin{array}{cccc}
\widetilde{G}_{j1} & & O
\\
 & \ddots &
\\
O & & \widetilde{G}_{jr}
\end{array}
\right] ~ (j=1,\ldots,m), \quad \widetilde{G}_{jk} \in \S^{t_{k}} ~ (k = 1,\ldots, r),
\end{eqnarray*}
where $s_{k}, ~ t_{k} \in \N ~ (k = 1,\ldots, r)$.
\par
To show this reformulation, we consider converting the decision variable matrix $R \in \R^{n \times r}$ of (NSDP$_{r}$) into a vector. Assume that $R$ can be written as $R = [u_{1} \cdots u_{r}]$, where $u_{k} \in \R^{n} ~ (k = 1,\ldots,r)$. For any $M \in \R^{n \times n}$, we have
\begin{eqnarray}
\langle M, RR^{\top} \rangle 
&=& {\rm tr}(MRR^{\top}) \nonumber
\\
&=& {\rm tr} \left( \left[
\begin{array}{c}
u_{1}^{\top}
\\
\vdots
\\
u_{r}^{\top}
\end{array}
\right] \left[
\begin{array}{ccc}
Mu_{1} & \cdots & Mu_{r}
\end{array}
\right] \right) \nonumber
\\
&=& \sum_{j=1}^{r} u_{j}^{\top} M u_{j} \nonumber
\\
&=& \left[
\begin{array}{ccc}
u_{1}^{\top} & \cdots & u_{r}^{\top}
\end{array}
\right] \left[
\begin{array}{ccc}
M & & O
\\
 & \ddots &
\\
O & & M
\end{array}
\right] \left[
\begin{array}{c}
u_{1}
\\
\vdots
\\
u_{r}
\end{array}
\right]. \label{RR_reform}
\end{eqnarray}
By exploiting (\ref{RR_reform}), it can be verified that (NSDP$_{r}$) is equivalent to (QECQP) with $d = nr$,
\begin{eqnarray*}
&& H = \left[
\begin{array}{ccc}
C & & O
\\
 & \ddots &
\\
O & & C
\end{array}
\right], \quad G_{j} = \left[
\begin{array}{cccc}
A_{j} & & O
\\
 & \ddots &
\\
O & & A_{j}
\end{array}
\right] ~ (j=1,\ldots,m).
\end{eqnarray*}
In a similar way to the reformulation of $R \in \R^{n \times r}$, we consider converting the decision variable matrix $S \in \L^{n \times r}$ of (T-NSDP$_{r}$) into a vector. Let $v_{j} \in \R^{n-k+1} ~ (k=1,\ldots,r)$ be the column of $S$ excluding the upper diagonal elements, i.e.,
\begin{eqnarray*}
S = \left[
\begin{array}{ccccc}
\cellcolor[gray]{.8}& 0 & 0 & \cdots &0
\\
\cellcolor[gray]{.8}& \cellcolor[gray]{.8} & 0 & & \vdots
\\
\cellcolor[gray]{.8}& \cellcolor[gray]{.8} & \cellcolor[gray]{.8} & \ddots &
\\
\cellcolor[gray]{.8}& \cellcolor[gray]{.8} & \cellcolor[gray]{.8} &  & 0
\\
\cellcolor[gray]{.8}v_{1} & \cellcolor[gray]{.8}v_{2} & \cellcolor[gray]{.8}v_{3} &\cdots& \cellcolor[gray]{.8}v_{r}
\\
\cellcolor[gray]{.8}& \cellcolor[gray]{.8} & \cellcolor[gray]{.8} & &\cellcolor[gray]{.8}
\\
\cellcolor[gray]{.8}& \cellcolor[gray]{.8} & \cellcolor[gray]{.8} & &\cellcolor[gray]{.8}
\\
\cellcolor[gray]{.8}& \cellcolor[gray]{.8} & \cellcolor[gray]{.8} & &\cellcolor[gray]{.8}
\\
\cellcolor[gray]{.8}& \cellcolor[gray]{.8} & \cellcolor[gray]{.8} & &\cellcolor[gray]{.8}
\\
\cellcolor[gray]{.8}& \cellcolor[gray]{.8} & \cellcolor[gray]{.8} & &\cellcolor[gray]{.8}
\end{array}
\right] \in \L^{n \times r}.
\end{eqnarray*}
Then, we also see that (T-NSDP$_{r}$) can be reformulated as (QECQP) with $d = nr - \frac{r(r-1)}{2}$,
\begin{eqnarray*}
&& H = \left[
\begin{array}{ccc}
C_{1} & & O
\\
 & \ddots &
\\
O & & C_{r}
\end{array}
\right], \quad G_{j} = \left[
\begin{array}{cccc}
A_{j1} & & O
\\
 & \ddots &
\\
O & & A_{jr}
\end{array}
\right] ~ (j=1,\ldots,m),
\end{eqnarray*}
where $C_{k}, A_{1k}, \ldots, A_{mk} ~ (k=1,\ldots,r)$ are matrices obtained by removing the first $k-1$ rows and columns of $C, A_{1}, \ldots, A_{m}$, respectively.
\par
From now on, we discuss how to solve (QECQP). Here, we utilize a local convergent SQP method \cite[Algorithm~18.1]{NoWr06} for (QECQP). To give its formal statement, we define some notation. Let $g \colon \R^{d} \to \R^{m}$ and $\nabla g \colon \R^{d} \to \R^{d \times m}$ be defined as
\begin{eqnarray*}
g(x) := \left[
\begin{array}{c}
g_{1}(x)
\\
\vdots
\\
g_{m}(x)
\end{array}
\right], \quad \nabla g(x) := \left[
\begin{array}{ccc}
\nabla g_{1}(x) & \cdots & \nabla g_{m}(x)
\end{array}
\right].
\end{eqnarray*}
Moreover, let $\mathscr{L} \colon \R^{d} \times \R^{m} \to \R$ be the Lagrangian associated with (QECQP) defined as 
\begin{eqnarray*}
\mathscr{L}(x, \mu) := f(x) - \langle \mu, g(x) \rangle.
\end{eqnarray*}
For completeness, we present below this local convergent SQP algorithm.
\begin{algorithm} \label{algorithm_for_NM}
An SQP method for {\rm (QECEP)}
\
\begin{description}

\item[Step {\rm 0}:] Choose an initial point $(x_{0}, \mu_{0}) \in \R^{d} \times \R^{m}$. Set $k := 0$.

\item[Step {\rm 1}:] Obtain the solution $(\xi^{\ast}, \zeta^{\ast})$ by solving
\begin{eqnarray*} \label{solution_subpro}
\left[
\begin{array}{cc}
\nabla_{xx}^{2} \mathscr{L}(x_{k}, \mu_{k}) & -\nabla g(x_{k})
\\
\nabla g(x_{k})^{\top} & O
\end{array}
\right] \left[
\begin{array}{cc}
\xi^{\ast}
\\
\zeta^{\ast}
\end{array}
\right] = - \left[
\begin{array}{cc}
\nabla f(x_{k})
\\
g(x_{k})
\end{array}
\right].
\end{eqnarray*}

\item[Step {\rm 2}:] Set $x_{k+1} := x_{k} + \xi^{\ast}$ and $\mu_{k+1} := \zeta^{\ast}$.

\item[Step {\rm 3}:] Set $k \leftarrow k + 1$, and go back to Step~{\rm 1}.
\end{description}
\vspace{-3mm}
\end{algorithm}

\section{Numerical experiments}
In this section, we report some numerical experiments using Algorithm~\ref{algorithm_for_NM}. All the programs were implemented with MATLAB R2020a and ran on a machine with Intel Core i9-9900k 3.60GHz CPU and 16GB of RAM. We basically compare the performance of the following three algorithms.
\begin{description}
\item[{\bf SDPT3:}] It is a well-known MATLAB solver for SDP \cite{ToToTu99}. We solve (SDP) by utilizing SDPT3.

\item[{\bf Algorithm~NSDP$_{r}$:}] It is an algorithm composed of SDPT3 and Algorithm~\ref{algorithm_for_NM} for (NSDP$_{r}$), where $r = \lceil (\sqrt{8m+1} - 1)/2 \rceil$. After SDPT3 finds an initial point of Algorithm~\ref{algorithm_for_NM}, which satisfies
\begin{eqnarray} \label{initial_criterion}
\max \left\{ \frac{\Vert {\cal A}(X_{0}) - b \Vert}{1 + \Vert b \Vert}, ~ \frac{\Vert C - {\cal A}^{\ast}(y_{0}) + Z_{0} \Vert_{{\rm F}}}{1 + \Vert C \Vert_{{\rm F}}} \right\} < 10^{-3},
\end{eqnarray}
Algorithm~\ref{algorithm_for_NM} solves (NSDP$_{r}$) by using $(X_{0},y_{0},Z_{0})$.

\item[{\bf Algorithm~T-NSDP$_{r}$:}] It is an algorithm composed of SDPT3 and Algorithm~\ref{algorithm_for_NM} for (T-NSDP$_{r}$), where $r = \lceil (\sqrt{8m+1} - 1)/2 \rceil$. After SDPT3 finds an initial point of Algorithm~\ref{algorithm_for_NM}, which satisfies condition (\ref{initial_criterion}), Algorithm~\ref{algorithm_for_NM} solves (T-NSDP$_{r}$) by using $(X_{0},y_{0},Z_{0})$.
\end{description}
\par
Moreover, the experiments mainly compare the computational time and the accuracy of solutions obtained by the above three algorithms. All test problems used here can be obtained in SDPT3 package. Throughout the experiments, we adopted the following as stopping conditions of Algorithm~\ref{algorithm_for_NM}:
\begin{eqnarray*} \label{termination_criterion}
E(X,y,Z) := \max \left\{ \frac{\Vert {\cal A}(X) - b \Vert}{1 + \Vert b \Vert}, ~ \frac{| \langle X, Z \rangle |}{1 + | \langle C, X \rangle | + | \langle b, y \rangle |} \right\} \leq \varepsilon ~~ {\rm or} ~~ k = 100,
\end{eqnarray*}
where we set $\varepsilon := 10^{-8}$ when comparing the computational time, and $\varepsilon := 0$ when comparing the accuracy of solutions. Similarly, concerning the stopping criteria of SDPT3, we set \texttt{gaptol} $= 10^{-8}$ and $\texttt{gaptol} = 0$, respectively, for the comparison of computational time.
\par
To check the accuracy of the solutions, the following indicators are utilized:
\begin{itemize}
\item the stopping criterion using $E(X,y,Z)$,
\item the infeasibility calculated by $\Vert {\cal A}(X) - b \Vert + \max \{ -\lambda_{\min}(C - {\cal A}^{\ast}(y)) ,0 \}$,
\item the duality gap given by $\langle C, X \rangle - \langle b, y \rangle$.
\end{itemize}

\noindent
{\bf The standard SDP}
\\
First of all, we solve the following standard SDP:
\begin{eqnarray} \label{standardSDP}
\begin{array}{ll}
\displaystyle \Mini_{X \in \S^{n}} & \langle C, X \rangle
\\
\subj & \langle A_{j}, X \rangle = b_{j} ~ (j=1,\ldots,m), ~ X \succeq O,
\end{array}
\end{eqnarray}
where $C$, $A_{j}~ (j=1,\ldots,m)$, and $b_{j} ~ (j=1,\ldots,m)$ were generated by \texttt{sdprand}, which is a command included in the package of SDPT3 that generates a random SDP.
\par
\begin{figure}[h]
\centering
\includegraphics[scale=0.7]{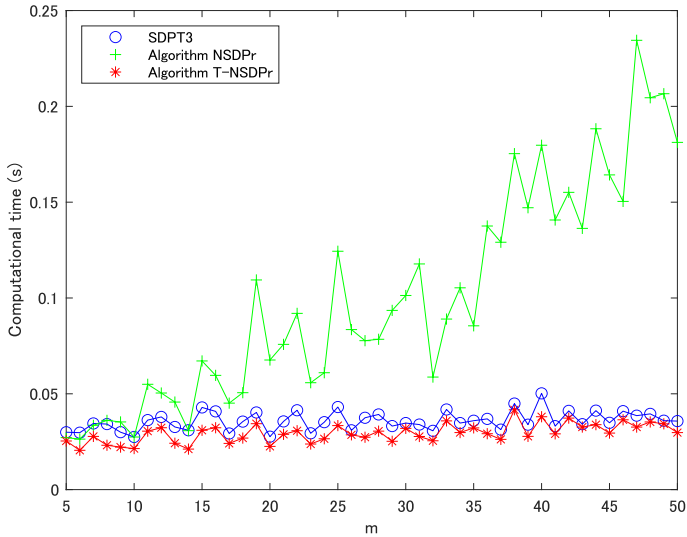}
\caption{Computational time obtained by each algorithm} \label{Exp1_total_time}
\end{figure}

In the first experiment, the dimension $n$ was set to $10$, and the number $m$ was incremented by $1$ from $5$ to $50$, and $5$ random problems were solved for each $m$. Figure~\ref{Exp1_total_time} indicates the average computational time of SDPT3, Algorithm~NSDP$_{r}$, and Algorithm~T-NSDP$_{r}$ for each $m$. The horizontal axis displays the number of constraints, that is $m$, and the vertical axis displays the computational time (in seconds). From Figure~\ref{Exp1_total_time}, we can confirm that Algorithm~T-NSDP$_{r}$ was faster than SDPT3 and Algorithm~NSDP$_{r}$ in most of the test problems. Moreover, even when $m$ was increased, the computational time of SDPT3 and Algorithm~T-NSDP$_{r}$ practically did not vary differently to Algorithm~NSDP$_{r}$. 
\par
In the second experiment, we set $n = 10$ and $m = 30$, and solved $10$ random problems. Table~\ref{prob1} illustrates the averages of the three indicators defined previously, i.e., the stopping criterion, the infeasibility, and the duality gap, obtained by the three algorithms. Clearly, by seeing all the three indicators, Algorithm~T-NSDP$_{r}$ was able to solve all the problems most accurately.
\begin{table}[h]
\begin{center}
\caption{The averages of the three indicators obtained by each algorithm} \label{prob1}
\begin{tabular}{cccccccccccc} 
\hline
 & SDPT3 & NSDP$_{r}$ & T-NSDP$_{r}$    
\\ \hline\hline
Stopping criterion & 7.4e-12	& 7.4e-12 & 6.3e-16	
\\ \hline
Infeasibility          & 4.8e-09	& 4.8e-09 & 6.1e-13
\\ \hline
Duality gap          & 6.8e-09	& 6.8e-09 & -3.8e-14
\\ \hline
\end{tabular}
\end{center}
\end{table}

\noindent
{\bf The max-cut problem}
\\
Let us now consider the max-cut problem below:
\begin{eqnarray*}
\begin{array}{ll}
\displaystyle \Mini_{X \in \S^{n}} & \displaystyle  \frac{1}{4}\langle B - {\rm diag}(Be), X \rangle
\\
\subj & [X]_{jj} = 1 ~ (j=1,\ldots,n), ~ X \succeq O,
\end{array}
\end{eqnarray*}
where $B$ is a weighted adjacency matrix of a graph generated by the command \texttt{graph} also included in SDPT3 package that generates random weighted adjacency matrices.
\par
\begin{figure}[h]
\centering
\includegraphics[scale=0.7]{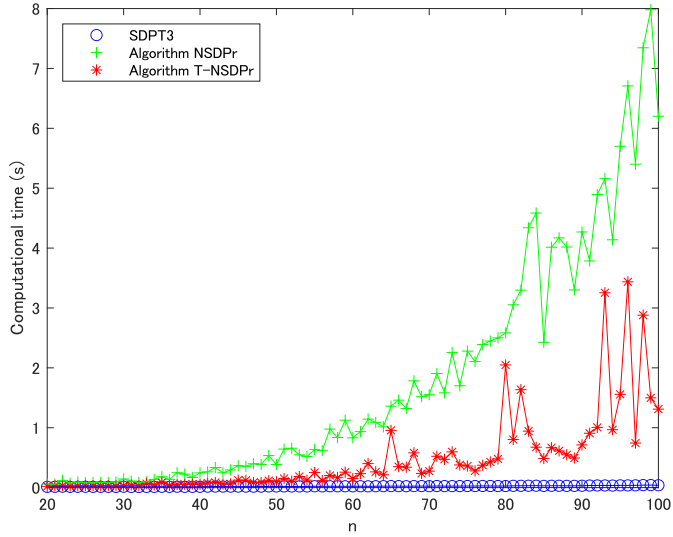}
\caption{Computational time obtained by each algorithm} \label{Exp2_total_time}
\end{figure}

We conducted the experiment related to the computational time as follows: The dimension $n$ was incremented one by one from $20$ up to $100$, and $5$ random problems were solved for each $n$. Figure~\ref{Exp2_total_time} shows the average computational time of SDPT3, Algorithm~NSDP$_{r}$, and Algorithm~T-NSDP$_{r}$ for each $n$. The horizontal axis represents the dimension of the decision variable $X$ and the vertical axis represents the computational time (in seconds). Figure~\ref{Exp2_total_time} shows that SDPT3 was the fastest algorithm. However, when the number of constraints is relatively small (less than $60$), Algorithm~T-NSDP$_{r}$ was competitive to SDPT3.
\par
We explain the experiment regarding the accuracy of solutions. The dimension $n$ was set to $50$, and $10$ random problems were solved. Table~\ref{prob2} gives the stopping criterion, the infeasibility, and the duality gap obtained by each algorithm. Although the results of Algorithm~T-NSDP$_{r}$ were inferior to those of Algorithm~NSDP$_{r}$ for many test problems, they were more accurate than those of SDPT3.
\begin{table}[h]
\begin{center}
\caption{The averages of three indicators obtained by each algorithm} \label{prob2}
\begin{tabular}{cccccccccccc} 
\hline
 & SDPT3 & NSDP$_{r}$ & T-NSDP$_{r}$    
\\ \hline\hline
Stopping criterion & 1.1e-13	& 1.2e-14 & 3.4e-14	
\\ \hline
Infeasibility          & 6.2e-12	& 5.1e-13 & 1.4e-12
\\ \hline
Duality gap          & 5.0e-11	& -1.3e-12 & 9.4e-13
\\ \hline
\end{tabular}
\end{center}
\end{table}

\noindent
{\bf The norm-minimization problem}
\\
Finally, we deal with the following norm-minimization problem:
\begin{eqnarray*}
\begin{array}{ll}
\displaystyle \Mini_{z \in \C^{m}} & \displaystyle \left\Vert \sum_{k=1}^{m} [z]_{k} B_{k} + B_{0} \right\Vert_{2},
\end{array}
\end{eqnarray*}
where $\mathbb{C}^{m}$ indicates the $m$-dimensional complex vector space, $B_{0}, B_{1}, \ldots, B_{m} \in \R^{p \times q}$ are constant matrices, which were generated by the MATLAB command \texttt{rand}. This problem can be reformulated as follows:
\begin{eqnarray*}
\begin{array}{cl}
\displaystyle \Mini_{t, x, y} & t,
\\
\subj & \displaystyle \sum_{k=1}^{m} [x]_{k} \left[
\begin{array}{cc}
O & B_{k}
\\
B_{k}^{\ast} & O
\end{array}
\right] + \sum_{k=1}^{m} [y]_{k} \left[
\begin{array}{cc}
O & iB_{k}
\\
(iB_{k})^{\ast} & O
\end{array}
\right] - t I \preceq - \left[
\begin{array}{cc}
O & B_{0}
\\
B_{0}^{\ast} & O
\end{array}
\right],
\end{array}
\end{eqnarray*}
where $i$ represents the imaginary unit, and the superscript $\ast$ denotes the conjugate transposition of a matrix. In this experiment, we solve the above reformulated problem.
\par
\begin{figure}[h]
\centering
\includegraphics[scale=0.7]{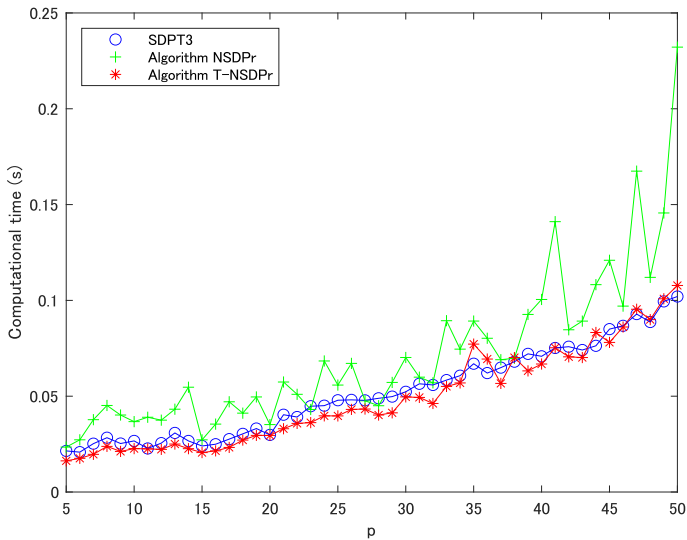}
\caption{Computational time obtained by each algorithm} \label{Exp3_total_time}
\end{figure}

In the experiment of the computational time, we set $m=10$ and $p=q$, and increased $p$ one by one from $5$ to $50$. In Figure~\ref{Exp3_total_time}, we provide the average computational time of SDPT3, Algorithm~NSDP$_{r}$, and Algorithm~T-NSDP$_{r}$ for each $p$. The horizontal axis represents the number of rows of constant matrices $B_{j} ~ (j=1,\ldots,m)$ and the vertical axis indicates the computational time (in seconds). Moreover, we solved $5$ random problems for each $p$. Figure~\ref{Exp3_total_time} shows that Algorithm~T-NSDP$_{r}$ was competitive to SDPT3. In particular, in the case where the dimension $p$ is less than $35$, Algorithm~T-NSDP$_{r}$ was superior to SDPT3 for most of the test problems. 

\par
We report results of the experiment associated with the accuracy of solutions. We set $p=q=50$ and $m=10$, and solved $10$ random problems. Table~\ref{prob3} shows the average accuracies of each algorithm. For all the test problems, Algorithm~T-NSDP$_{r}$ could solve them accurately compared with SDPT3 and Algorithm~NSDP$_{r}$.
\begin{table}[h]
\begin{center}
\caption{The average of the three indicators obtained by each algorithm} \label{prob3}
\begin{tabular}{cccccccccccc} 
\hline
 & SDPT3 & NSDP$_{r}$ & T-NSDP$_{r}$    
\\ \hline\hline
Stopping criterion & 6.3e-14	& 6.3e-14 & 3.6e-15	
\\ \hline
Infeasibility          & 2.7e-13	& 2.7e-13 & 1.9e-14
\\ \hline
Duality gap          & 4.2e-13	& 4.2e-13 & 7.2e-15
\\ \hline
\end{tabular}
\end{center}
\end{table}

\section{Conclusion}
In this paper, we have proposed a new nonlinear optimization model (T-NSDP$_r$) for (SDP), which can overcome the drawbacks of the existing model (NSDP$_{r}$) presented by Burer and Monteiro~\cite{BuMo03,BuMo05}. Since the decision variable spaces of (NSDP$_{r}$) and (T-NSDP$_{r}$) are respectively $\R^{n \times r}$ and $\L^{n \times r} = \{ S \in \R^{n \times r} \colon [S]_{ij} = 0 ~ \mbox{if} ~ i<j \}$, the proposed model is less than the existing one by $\frac{r(r-1)}{2}$-dimensions. Moreover, this dimensional reduction produce a beneficial result that (T-NSDP$_{r}$) has a strict local optimum under some appropriate conditions, whereas (NSDP$_{r}$) has no strict local one. Hence, we can expect that second-order methods for (T-NSDP$_{r}$) can quickly obtain solutions. Furthermore, we have conducted some numerical experiments which demonstrate that solving (T-NSDP$_{r}$) with an SQP method can efficiently find a more accurate solution.
\par
A future work is to provide sufficient conditions under which second-order sufficient conditions of (T-NSDP$_{r}$) hold.


\begin{thebibliography}{}
\bibitem{Ba95}
Barvinok, A.: Problems of distance geometry and convex properties of quadratic maps. Disc. Comput. Geom. {\bf 13}, 189–202 (1995)

\bibitem{Be09}
Bernstein, D.S.: Matrix Mathematics: Theory, Facts, and Formulas, 2nd edn, Princeton University Press, Princeton (2009)

\bibitem{BuMo03}
Burer, S., Monteiro, R.D.C.: A nonlinear programming algorithm for solving semidefinite programs via low-rank factorization. Math. Program. Ser. B {\bf 95}, 329--357 (2003)

\bibitem{BuMo05}
Burer, S., Monteiro, R.D.C.: Local minima and convergence in low-rank semidefinite programming. Math. Program. Ser. A {\bf 103}, 427--444 (2005)

\bibitem{MoOrSv14}
Monteiro, R.D.C., Ortiz, C., Svaiter, B.F.: Implementation of a block-decomposition algorithm for solving large-scale conic semidefinite programming problems. Comput. Optim. Appl. {\bf 57(1)}, 45--69 (2014)

\bibitem{NoWr06}
Nocedal, J., Wright, S.J.: Numerical Optimization, Springer, Berlin (2006)

\bibitem{Pa98}
Pataki, G.: On the rank of extreme matrices in semidefinite programs and the multiplicity of optimal eigenvalues. Math. Oper. Res. {\bf 23}, 339–358 (1998)

\bibitem{PoReWi06}
Povh, J., Rendl, F., Wiegele, A.: A boundary point method to solve semidefinite programs. Computing {\bf 78(3)}, 277--286 (2006)

\bibitem{St99}
Sturm, J.F.: Using SeDuMi 1.02, a MATLAB toolbox for optimization over symmetric cones. Optim. Methods Softw. {\bf 11(1-4)}, 625--653 (1999)

\bibitem{SuToYa16}
Sun, D., Toh, K.C., Yang, L.: An efficient inexact ABCD method for least squares semidefinite programming. SIAM J. Optim. {\bf 26(2)}, 1072--1100 (2016)

\bibitem{To01}
Todd, M.J.: Semidefinite optimization. Acta Numer. {\bf 10}, 515--560 (2001)

\bibitem{ToToTu99}
Toh, K.C., Todd, M.J., T\"ut\"unc\"u, R.H.: SDPT3  -- a MATLAB software package for semidefinite programming, version 1.3. Optim. Methods Softw. {\bf 11(1-4)}, 545--581 (1999)

\bibitem{Va96}
Vandenberghe, L., Boyd, S.: Semidefinite programming. SIAM Rev. {\bf 38(1)}, 49--95 (1996)

\bibitem{WeGoYi09}
Wen, Z., Goldfarb, D., Yin, W.: Alternating direction augmented Lagrangian methods for semidefinite programming. Math. Program. Comput. {\bf 2}, 203--230 (2010)

\bibitem{WoSaVa00}
Wolkowicz, H., Saigal, R., Vandenberghe, L. (eds.): Hnadbook of Semidefinite Programming: Theory, Algorithms and Applications, Kluwer International Series in Operations Research and Management Science. Kluwer, Boston (2000)

\bibitem{YaFuKo03}
Yamashita, M., Fujisawa, K., Kojima, M.: Implementation and evaluation of SDPA 6.0 (semidefinite programming algorithm 6.0). Optim. Methods Softw. {\bf 18(4)}, 491--505 (2003)

\bibitem{YaSuTo15}
Yang, L., Sun, D., Toh, K.C.: SDPNAL+: a majorized semismooth Newton-CG augmented Lagrangian method for semidefinite programming with nonnegative constraints. Math. Program. Comput. {\bf 7}, 331--366 (2015)

\bibitem{ZhSuTo10}
Zhao, X.Y., Sun, D., Toh, K.C.: A Newton-CG augmented Lagrangian method for semidefinite programming. SIAM J. Optim. {\bf 20(4)}, 1737--1765 (2010)
\end{thebibliography}
\end{document}